\documentstyle[11pt]{article}

\input amssym.def 
\input amssym

\newcommand{\be}{\begin{eqnarray*}}
\newcommand{\ee}{\end{eqnarray*}}

\newcommand{\equ}{\equiv 0 (mod \, D)}
\newcommand{\ch}{( D, \, H )}
\newcommand{\mr}{H^{\#}}

\newcommand{\mm}{H}

\newcommand{\D}{U\times \f{h}}

\newcommand{\al}{\alpha}

\newcommand{\lam}{\lambda}
\newcommand{\va}{\varphi}
\newcommand{\rr}{\Lambda^{\#}}
\newcommand{\rrr}{[\Lambda , \Lambda ] ^{\#}}
\newcommand{\rt}{\Lambda}

\newcommand{\adh}{$ad_{{\frak h}}$ - invariant}
\newcommand{\f}[1]{{\frak #1}}  
\newcommand{\plam}{{\partial  \over \partial \lambda_i}}
\newcommand{\plamm}[1]{{\partial {#1} \over \partial \lambda_i}}

\newcommand{\4}[4]{(#1,\,#2 \,;\,#3,\,#4)}

\newcommand{\lon}{\longrightarrow}
\newcommand{\complex}{{\Bbb C}}

\newcommand{\half}{\frac{1}{2}}

\newcommand{\cinf}{C^{\infty}}

\newtheorem{thm}{Theorem}[section]
\newtheorem{lem}[thm]{Lemma}
\newtheorem{cor}[thm]{Corollary}
\newtheorem{pro}[thm]{Proposition}

\newcommand{\pf}{\noindent{\bf Proof.}\ }
\newcommand{\qed}{\begin{flushright} $\Box$\ \ \ \ \  \end{flushright}}

\newcommand{\rsd}
{ {\frak g}= {\frak h} \oplus \sum_{\alpha \in \Delta}
{\frak g_{\alpha}}}
\newcommand{\frakg}{{\frak g}}
\newcommand{\frakh}{{\frak h}}
\newcommand{\frakk}{{\frak k}}

\newcommand{\gm}{\Gamma }
\newcommand{\T}{\tau}

\newcommand{\epe}{\mbox{$\epsilon$}}
\newcommand{\ea}{\mbox{$E_{\alpha}$}}
\newcommand{\eb}{\mbox{$E_{-\alpha}$}}

\begin{document}

\title{{\bf Dirac structures and dynamical r-matrices}}
\author{Zhang-Ju Liu 
\thanks{Research partially supported
by NSF of China and the Research Project of ``Nonlinear Science".}\\ 
	Department  of Mathematics \\
	Peking  University \\
	Beijing, 100871, China\\
        {\sf email: liuzj@pku.edu.cn}\\
        Ping Xu 
         \thanks{ Research partially supported by NSF
       grant DMS97-04391. }\\ 
        Department of Mathematics\\
         Pennsylvania State University \\
         University Park, PA 16802, USA \\
	{\sf email: ping@math.psu.edu }}

\date{March, 1999}
\maketitle
\begin{abstract}
The purpose of this paper is to establish a connection between
various subjects such as dynamical $r$-matrices, Lie bialgebroids,
and Lagrangian subalgebras. Our  method relies 
on the theory of Dirac structures developed in \cite{LWX}
\cite{LWX1}.  In particular, we give a new method of  
classifying  dynamical $r$-matrices of simple Lie algebras
$\frakg$, and prove that dynamical $r$-matrices are in one-one correspondence
with certain Lagrangian subalgebras of $\frakg \oplus \frakg $.
\end{abstract}

\section{Introduction}

Recently, there has been a great  deal of interest in  the so called 
{\em Classical Dynamical Yang-Baxter Equation} (here after 
{\em CDYBE}):

\begin{equation}
\label{cdybe}
 Alt(dr) \, + \,
[r^{12},r^{13}] + [r^{12}, r^{23}] + [r^{13}, r^{23}] \, = \, 0,
\end{equation}
where  $r (\lambda ): \frakh^{*} \rightarrow \f{g} \otimes \f{g}$ is
a meromorphic function, and 
$\f{g}$ is  a complex simple Lie algebra with Cartan
subalgebra $\f{h}$.
 When $r$ is a constant function, Equation
(\ref{cdybe}) reduces to the usual classical Yang-Baxter equation,
and therefore a classical $r$-matrix is a special solution.
Assume that   $r$ is a solution,  and
 that $r+r^{21}=\epsilon \Omega$, where
$\Omega\in (S^{2}\frakg)^{\frakg}$ is the  Casimir element
corresponding to the Killing form,
and $\epsilon$ is a constant usually called the coupling constant.
Then the skew-symmetric part of $r$ satisfies 
  the following modified  CDYBE:
\begin{equation}
Alt(dr) + \half [ r, r] =  {\epe^2 \over 4}\,
[\Omega^{12},\Omega^{23}]\,  \in (\wedge^{3}\frakg )^{\frakg},
 \label{rr}
\end{equation}
where $ [ \cdot, \cdot ]$ is the Schouten bracket on $ \wedge^{*}\f{g}$.

In this paper,  by a {\it dynamical $r$-matrix}, we mean
a  meomorphic function  $r:\frakh^{*} \rightarrow \f{g} \wedge \f{g}$
satisfying:
\begin{enumerate}
\item $[h, r (\lambda )] = 0, \ \ \ \forall h \in \f{h}$, and
\item  $r$ satisfies the modified  CDYBE (\ref{rr}).
\end{enumerate}
The first assumption is often referred to as the zero
weight condition \cite{EV}. Here we are mainly interested in dynamical
$r$-matrix with nonzero coupling constant. In this case,
by multiplying by a constant, we may always assume
that $\epsilon=2$. In the sequel, we will 
always make this assumption when referring to a dynamical
$r$-matrix unless otherwise  specified.
 
Classical dynamical $r$-matrices have appeared in various  contexts
in mathematical physics,
for instance,  in Knizhnik-Zamolodchikov-Bernard equation \cite{F}, and
in the study of integrable systems such as Caloger-Moser
systems \cite{Avan} \cite{BAB} \cite{BAB1}.   A
 classification of dynamical r-matrices
for simple Lie algebras  was obtained  by Etingof and Varchenko
in \cite{EV}.  An example of such a dynamical $r$-matrix is 
\[
r(\lambda) \, = 
 \sum _{\alpha \in \Delta_+} \,
\coth ( < \alpha, \lambda >) \ea \wedge  \eb,
\]
where $\Delta_{+}$ is the set of positive
roots of $\frakg$ with respect to $\f{h}$, the
$\ea$ and $\eb$'s are root vectors, and
$\coth (x)  =  {e^x + e^{-x} \over e^x - e^{-x}}$
is the hyperbolic cotangent function.
Moreover, it is proved that in \cite{EV}
 dynamical $r$-matrices correspond
to Poisson groupoids just  as
 classical r-matrices integrate  to Poisson 
groups in Drinfeld theory \cite{LW} \cite{S}. 
The  corresponding Lie bialgebroids, as the
 infinitesimal invariants, were studied by
Bangoura and Kosmann-Schwarzbach \cite{BK-S}.

It is well known that
there are many  ways of producing
 a classical $r$-matrix. A  natural method is
 via  Lie bialgebras using    Manin triples.
 For instance,
for the Lie bialgebra of  the
standard $r$-matrix $r_{0}=\sum_{\alpha \in \Delta_+}\ea \wedge  \eb$,
the corresponding Manin triple is  
($\frakg \oplus \frakg , \ \frakg_1, \ \frakg_{2})$, where
$\frakg_1 \subset \frakg$ is  the diagonal while $\frakg_2$
 is the subalgebra $\{(h+X_{+}, -h+ X_{-})|h\in \frakh , X_{\pm}
\in n_{\pm}\}$. Here $n_{\pm} \subset \frakg$ are maximal nilpotent
 subalgebras.
It  is thus  natural to ask

\begin{quote}
{\bf Problem 1}. Does  there exist such  an analogue
 for dynamical r-matrices?
In particular, what is the  double of the Lie bialgebroid 
corresponding to a  dynamical r-matrix?
\end{quote}

Recently, Lu has found an interesting connection
between  dynamical $r$-matrices and  Poisson homogeneous spaces \cite{Lu}.
More precisely, Lu showed that a dynamical $r$-matrix
gives rise to a family of Poisson homogeneous $G$-spaces  $G/H$
parameterized by $\lambda $,
 where $G$ is the Poisson group defined by the standard classical
$r$-matrix $r_{0}$ with the same coupling 
constant (i.e., constant solution of Equation
(\ref{rr})), and $H$ is the subgroup of $G$ having
Lie algebra $\frakh$.  Clearly, the Poisson homogeneous spaces
corresponding to different $\lambda$, 
must  be related in some way that
 should  reflect the dynamical property of the dynamical $r$-matrix.
This leads to our

\begin{quote}
{\bf Problem 2}.  Given a family of Poisson homogeneous $G$-spaces  $G/H$
parameterized by $\lambda$,  what criteria
 will guarantee that it arises from a dynamical r-matrix?
\end{quote}

The infinitesimal  object of the Poisson group $G$ is
the Lie bialgebra $(\frakg , \frakg^{*}, r_0 )$ 
generated by the classical $r$-matrix $r_0$. According
to Drinfeld \cite{dr:poisson}, Poisson homogeneous $G$-spaces
are in one-one correspondence with  Lagrangian subalgebras
of the double Lie algebra $\f{d}$, which is
isomorphic to the direct sum Lie algebra $\frakg \oplus \frakg$. So
an equivalent formulation of Problem 2 is 

\begin{quote}
{\bf Problem 3}.  Let $W(\lambda )\subset \f{d}$ be
a family of Lagrangian subalgebras. When will  this family
of Lagrangian subalgebras be induced from a dynamical r-matrix?
\end{quote}

In fact Lu showed that these   Poisson homogeneous $G$-spaces
include all the Poisson homogeneous $G$-spaces of the form $G/H$.
 This suggests that dynamical $r$-matrices and  Lagrangian subalgebras
of $\frakg \oplus \frakg$  should  be  intrinsically
related in some manner.  On the other hand, a general
  classification of Lagrangian
 subalgebras of $\f{d}$ has been  obtained  by Karolinsky \cite{Ka}, which
does not seem to have  an obvious connection with the work
of Etingof and Varchenko \cite{EV}. Therefore
it is natural to ask

\begin{quote}
{\bf Problem 4}.  What is the precise relation 
between dynamical r-matrices and Lagrangian subalgebras
of $\frakg \oplus \frakg $?
\end{quote}

Th purpose  of this paper is to understand the
intrinsic connection between various objects such as dynamical
$r$-matrices, Lagrangian subalgebras, and Lie bialgebroids (see \cite{W}).
In particular, our work is 
 motivated by the above  questions. Our idea is to  use
Dirac structure theory developed in \cite{LWX} \cite{LWX1}.
The starting point is a  simple Courant algebroid (see Section 3):
$(TU\oplus T^* U)\times (\frakg \oplus \frakg )$,
which can be considered as an analogue  of the direct
sum Lie algebra $\frakg \oplus \frakg $ in the
algebroid context, where $U\subset \frakg^*$ 
is an open subset. We   analyze a class of
 Dirac structures of this Courant algebroid
which  are induced from dynamical $r$-matrices. 
This study leads to a new method of classification
of dynamical $r$-matrices for simple Lie algebras.
One advantage of our approach is that the Cayley
transformation, which turns out to be important
in classification theory \cite{Sch}, appears quite naturally.
We hope  that our method may shed  new light on the 
classification scheme of more general dynamical $r$-matrices \cite{Avan},
and that of dynamical $r$-matrices for compact Lie algebras.
This  discussion is the   main  topic  of  Section 4. 
In Section 5, we  show that Lagrangian subalgebras of
$\frakg \oplus \frakg$ whose intersection with the diagonal  are
equal to $\frakh$,  are in one-one correspondence with 
dynamical $r$-matrices with zero gauge term.
This relates  the results  of Karolinsky and Lu with that
of Etingof and Varchenko in an explicit way.
 Moreover, we prove that given a point
 $\mu \in \frakg^*$,  any such  Lagrangian subalgebra $W_0$  admits a unique
extension to a family of  Lagrangian subalgebras  $W(\lambda )$ with
$W(\mu ) =W_0$,  governed by a dynamical $r$-matrix. 
In a certain sense, this is similar to an initial value problem
of a first order o.d.e.
Section 2 contains some basic facts concerning
Lie bialgebroids and Courant algebroids. And Section 3
is devoted to the discussion on the connection between
dynamical $r$-matrices and Lie bialgebroids.

{\bf  Acknowledgments.}  In addition to the funding sources mentioned
in the first footnote, we would like to thank several institutions
for their hospitality while work on this project was being done:
IHES,  and Peking University (Xu); Penn State University (Liu).  Thanks go
also to Yvette Kosmann-Schwarzbach, and Jiang-hua Lu
for their helpful comments. Especially, we are grateful
to Lu for allowing us  to have access to her unpublished
manuscript \cite{Lu}.

\section{Preliminaries}

In this section, we recall some basic facts concerning Lie bialgebroids
and Dirac structures.

 A Lie bialgebroid is a pair of Lie algebroids ($A$, $A^*$) satisfying
 the following compatibility condition (see \cite{MX94} and
 \cite{K-S95}):
\begin{equation}
\label{eq:1}
d_{*}[X, Y]=[d_{*}X, Y]+[X, d_{*}Y], \ \ \forall X ,Y \in \Gamma (A),
\end{equation}
where the differential $d_*$ on $\Gamma(\wedge^*A)$
comes from the Lie algebroid structure on $A^*$.

Given  a Lie algebroid $A$  over $P$ with anchor $a$ and a section
$\Lambda \in \Gamma (\wedge^2 A )$,
 Denote by $\Lambda^\#$ the bundle map $A^{*}\lon A$ defined by
$\rr (\xi) (\eta )=\Lambda  (\xi , \eta ), \forall \xi , \eta \in
 \Gamma (A^{*})$.  Introduce  a bracket on
$\Gamma (A^* )$ by 
\begin{equation}
\label{eq*}
[\xi , \eta ]_{\Lambda} = L_{\rr \xi}\eta -L_{\rr \eta}\xi -d[\Lambda  (\xi , \eta )].
\end{equation}

By $a_*$ we denote the composition  $a\circ \rr :A^*  \lon TP$.

\begin{thm}
\label{LX}
$A^*$ with the bracket   and   anchor $a_*$  above
becomes   a Lie algebroid iff
\begin{equation}
\label{eq:2}
L_{X}[\rt, \rt ]= [ X, [\rt, \rt]] = 0, \  \ \  \forall X \in \Gamma (A).
\end{equation}
\end{thm}
\pf In \cite{LX96}, we proved this result with one more condition:
$a \circ \rrr =0$, which is equivalent to $[ f, [\rt,\rt]] = 0, \, \,
\forall f \in C^{\infty}(P)$. But in fact this last condition 
is a consequence of
Equation (\ref{eq:2}). To see this,  by replacing $X$ with $fX$
in Equation (\ref{eq:2}), one obtains
$[fX, [\Lambda , \Lambda ]]=0$. It thus follows that
$X\wedge [f,  [\Lambda ,\Lambda ]]=0$, $\forall X\in \Gamma (A)$,
which implies that  $[f, [\Lambda , \Lambda ]]=0$.
\qed
In this case, the induced  differential  $d_{*}: \Gamma (A )
\lon \Gamma (\wedge^{2} A)$ is simply given by $d_{*} X=[\Lambda , X ]$, 
$\forall X\in \Gamma (A )$.
Thus the   compatibility condition, Equation  (\ref{eq:1}),
  is satisfied automatically.  So $(A, A^* )$ is a Lie bialgebroid,
 called {\em  coboundary  Lie bialgebroid}.  $\Lambda$
is also called an $r$-matrix by abuse of notations.
 When $P$ reduces to a point, i.e., $A$ is a Lie
algebra,  Equation (\ref{eq:2}) is equivalent to that $[\rt , \rt ]$ is
$ad$-invariant, i.e, $\rt$ is a classical  $r$-matrix in
the ordinary sense.
On the other hand, when $A$ is the tangent
bundle $TP$ with the  standard Lie algebroid structure,
Equation (\ref{eq:2})  is  equivalent
to that $[\rt , \rt ]=0$, i.e., $\rt$ is a Poisson tensor.

Given a Lie bialgebroid ($A$,  $A^{*}$)  over the base
$P$,  with anchors $a$ and $a_{*}$ respectively,
let $E$ denote their   vector bundle direct sum:
$E=A\oplus A^{*}$.
On   $E$, there exists  a   natural non-degenerate symmetric
bilinear form:

\begin{equation}
\label{eq:inner}
(X_{1}+\xi_{1} , X_{2}+\xi_{2}) = 
\half (\langle \xi_{1},  X_{2}  \rangle + \langle \xi_{2} ,  X_{1}\rangle ).
\end{equation}

In \cite{LWX},    we  introduced a bracket  on $\Gamma (E)$, called 
{\em Courant bracket}:

\begin{equation}
\label{eq:double}
\begin{array}{lll}
[e_{1}, e_{2}] &= &\{ [X_{1}, X_{2}]+L_{\xi_{1}}X_{2}-L_{\xi_{2}}X_{1}-\half 
d_{*}  (\langle \xi_{1},  X_{2}  \rangle - \langle \xi_{2} ,  X_{1}\rangle )\}\\
&& + \{ [\xi_{1} , \xi_{2}]+L_{X_{1}}\xi_{2}-L_{X_{2}}\xi_{1} + \half
d (\langle \xi_{1},  X_{2}  \rangle - \langle \xi_{2} ,  X_{1}\rangle ) \},
\end{array}
\end{equation}
where $e_{1}=X_{1}+\xi_{1}$ and $e_{2}=X_{2}+\xi_{2}$.
Let $\rho : E\lon TP$ be the bundle map
$\rho =a +a_{*}$. That is,
\begin{equation}
\rho (X+\xi )=a(X)+a_{*} (\xi  ) , \ \ \forall X\in \Gamma (A) \mbox{ and }
\xi \in \Gamma (A^{*}).
\end{equation}

For a Lie bialgebra $(\f{g} , \f{g}^{*})$,
the  bracket (\ref{eq:double}) reduces to the  well known Lie
bracket  on   the double $\f{g} \oplus \f{g}^{*}$.
On the other hand, if $A$ is the tangent bundle Lie algebroid
$TM$ and $A^{*}=T^{*}M $  with zero bracket,
then Equation (\ref{eq:double}) takes  the form:
\begin{equation}
\label{eq:courant}
[X_{1}+\xi_{1} , X_{2}+\xi_{2}]=[X_{1} ,X_{2}]
+\{ L_{X_{1}}\xi_{2} -L_{X_{2}}\xi_{1} + \half
d (\langle \xi_{1},  X_{2}  \rangle - \langle \xi_{2} ,  X_{1}\rangle ) \}.
\end{equation}
This is the bracket first introduced by Courant \cite{Courant:1990}.
In general, $E$ together with this bracket
and the  bundle map $\rho$  satisfies certain properties as outlined in
the following:

\begin{thm} \cite{LWX}
\label{thm:LWX}
Given a Lie bialgebroid $(A, A^{*})$, let $E=A\oplus A^*$.
Then $E$,  together with
the  non-degenerate  symmetric  bilinear form
 $( \cdot , \cdot )$,  the  skew-symmetric
bracket $[\cdot , \cdot ]$ on $\Gamma (E)$
and the bundle map $\rho :E\lon TP$ as introduced above,
satisfies the following    properties:
 \begin{enumerate}
\item For any $e_{1}, e_{2}, e_{3}\in \Gamma (E)$,
$[[e_{1}, e_{2}], e_{3}]+c.p.=\f{D}  T(e_{1}, e_{2}, e_{3});$
\item  for any $e_{1}, e_{2} \in \Gamma (E)$,
$\rho [e_{1}, e_{2}]=[\rho e_{1}, \rho  e_{2}];$
\item  for any $e_{1}, e_{2} \in \Gamma (E)$ and $f\in C^{\infty} (P)$,
$[e_{1}, fe_{2}]=f[e_{1}, e_{2}]+(\rho (e_{1})f)e_{2}-
(e_{1}, e_{2})\f{D} f ;$
\item $\rho \circ \f{D} =0$, i.e.,  for any $f, g\in C^{\infty}(P)$,
$(\f{D} f,  \f{D}  g)=0$;
\item for any $e, h_{1}, h_{2} \in \Gamma (E)$,
  $\rho (e) (h_{1}, h_{2})=([e , h_{1}]+\f{D} (e ,h_{1}) ,
h_{2})+(h_{1}, [e , h_{2}]+\f{D}  (e ,h_{2}) )$,
\end{enumerate}

where
\begin{equation}
\label{eq:T0}
 T(e_{1}, e_{2}, e_{3})={ 1 \over 3} ([e_{1}, e_{2} ], e_{3})+c.p.,
\end{equation}
and
$\f{D} :  C^{\infty}(P)\lon \Gamma (E)$
is  the map $\f{D}=d_{*}+d $.
\end{thm}

$E$ is called the     double of the Lie bialgebroid $(A, A^* )$.
In general, a vector bundle $E$ equipped with the
 above structures is  called  a  {\em  Courant algebroid}
 \cite{LWX}.

In this paper, we are mainly interested in a special
gauge Lie algebroid  $A=TM\times \f{g} $, where
$\frakg$ is a Lie algebra. Clearly
 $A$ is a Lie algebroid over $M$ with  anchor being   the
projection $ p: A \lon TM$. As for the bracket, note
that any section of $A$ can always be
written as the sum of a vector field  and a $\f{g}$-valued function
on $M$. The bracket of  such two sections is given by:
\begin{equation}
\label{gb}
  [X + \xi, Y +\eta] = [ X,Y] + [ \xi, \eta] + L_X \eta -L_Y \xi,
   \, \,\, X ,Y \in  \chi (M), \, \,  \xi, \eta \in C^{\infty}(M, \f{g}),
\end{equation}
where the bracket of two vector fields is the usual bracket and the bracket
$[ \xi, \eta]$ is the pointwise bracket.  



Let  $r\in \wedge^{2} \f{g}$, which can be considered 
  as a constant section
of $\wedge^2 A$. Then
\begin{pro}
\label{pro:gauge}
$( A, A^*, r)$ is a coboundary 
 Lie bialgebroid iff $ [r, r] $ is $ad$-invariant,
 i.e., iff $( \f{g}, \f{g^*}, r)$ is a coboundary Lie bialgebra.
\end{pro}

In this case, the bracket for sections of  
$A^* (\cong  T^*M \times \f{g^*} )$ is given by
\begin{equation}
\label{gb*}
  [ \alpha + \xi, \beta +\eta] =[\xi, \eta ],
   \, \,\, \alpha , \beta \in  \Omega^1 (M), \, \,  \forall
 \xi, \eta \in C^{\infty}(M, \f{g^*}),
\end{equation}
where the right hand  side bracket is pointwise bracket on $\f{g^*}$.
The corresponding  double  is the vector bundle
$$
E = A \oplus A^*  \cong (TM \oplus T^*M) \times ( \f{g} \oplus\f{g^*}),
$$
where  the Courant  bracket can be described quite simply.
 On the   subbundle $TM \oplus T^*M$, the bracket is
just Courant's original bracket: Equation (\ref{eq:courant}), while
  for two elements of the double Lie algebra $ \f{g} \oplus \f{g^*}$  considered as constant sections of
$E$, the bracket is pointwise bracket. One should however
 note that   the
subbundle $ M \times (\f{g} \oplus \f{g^*})$ is not  closed
under
  the Courant bracket (\ref{eq:double}),  since  the third  property
in Theorem \ref{thm:LWX} implies that
\begin{equation}
\label{eq:fg}
[fe_1 , ge_2 ] = (fdg - gdf) ( e_1, e_2 ) \, + fg [ e_1, e_2 ],
 \,  \, \forall f, g \in C^{\infty}(M),  \,  \forall e_1,e_2
  \in \f{g} \oplus \f{g^*},
\end{equation}
 where $fdg - gdf \in \Omega^1( M)$.
 On the other hand, for $ X + \alpha\in \Gamma (TM \oplus T^*M)$,
 $f \in  C^{\infty}(M) $ and $ e \in \f{g} \oplus \f{g^*}$,
 we have
\begin{equation}
\label{eq:fe}
[ X + \alpha , fe ] = L_X (fe) = (Xf) e.
\end{equation}
 These formulas will be needed later on in Section  4.

Given a Courant algebroid $E$, a  {\em  Dirac structure}    is
a subbundle $L\subset E$
which is maximally isotropic with respect to
the symmetric bilinear form $( \cdot , \cdot )$
and  is  integrable in the sense that
 $\Gamma (L)$ is closed under the bracket $[\cdot , \cdot ]$.
 There are two   important  classes of Dirac structures studied in
 \cite{LWX}. One is the Dirac structures induced by
Hamiltonian operators, and the other is the so called
null Dirac structures. Let us  briefly recall their definitions
below.

Let $ \mm \in \gm (\wedge^2 A)$ and denote $ \mr : A^* \lon A$  the
induced  bundle map. Then the  graph of $\mr$ ,
$$  \gm_{ \mm} =\{ \mr \xi + \, \xi | \, \forall \xi \in A^* \},$$
defines a maximal isotropic subbundle of $A \oplus A^*$.
$\gm_{\mm}$ is a Dirac subbundle iff $H$
satisfies the Maurer-Cartan type equation: 
\begin{equation}
\label{eq:M}
d_{*}H+\half [H, H]=0.
\end{equation}
In this case we call $H$  
a {\it Hamiltonian operator}. Another interesting  class of Dirac structures 
is the so  called {\em null} Dirac structures,
which can be  characterized as follows.  
Let  $D \subseteq A$ be a subbundle, and $ D^{\perp}
\subseteq A^*$ its conormal  subbundle.
Consider  $ L = D \oplus D^{\perp} \subset A \oplus A^* $.
 Then $L$ is a Dirac structure iff  $D$   and $D^{\perp}$ are
Lie subalgebroids of $A$ and $A^* $, respectively .
In this case $L$ is called a null Dirac structure.

A more general construction of Dirac structures is
via the so called characteristic pairs \cite{Liu}.
  Let $D \subseteq A$ be   a subbundle and 
  $H \in \Gamma (\wedge^2 A)$. Define
\begin{equation}
\label{eq:L}
  L = \{ \,X + \mr \xi + \, \xi \,| \,  \forall
 X \in D , \xi \in D^{\perp} \} = D \oplus graph (\mr |_{D^{\perp}}),
\end{equation}
where $D^{\perp} \subseteq A^*$  is  the conormal subbundle  of $D$.
Clearly, $L$ is  a maximal isotropic subbundle of $A \oplus A^*$.
The pair $(D, H)$ is called a {\em characteristic pair } of $L$.

Conversely,  any  maximal isotropic subbundle $L\subset A$ such that 
$L \cap A $ is of constant rank can always  be described
by  such a characteristic pair.
Note that two  characteristic pairs   $(D_1, H_1)$ and   $(D_2, H_2)$
  define the  same  subbundle  $L$ by Equation (\ref{eq:L}) iff
$$
D_1 = D_2, \ \ \ \
  \mbox{and}  \ \ \ \  pr ( H_1) = pr( H_2 ),
 \,  i.e., \, H_1 - H_2 \ \ 
 \equ ,$$
where $pr$ denotes the projection $A \lon A/D$ and
its induced map $\gm (\wedge^{*}A)\lon \gm (\wedge^{*}A/D)$.
 In the above equation as well as in  the sequel,
  a section $ \Omega \in \gm (\wedge^* A)$ is said
equal to zero module $D$, denoted   as $\Omega \equ$,
if its projection under $pr$ vanishes in  $\gm (\wedge^* A/D)$.
Even though $L$ is related only to $ pr(\mm)  \in \gm (\wedge^2 A/D)$
instead of \, $ \mm$ itself,
it is still  convenient to characterize the integrability conditions of $ L$
in terms of  $\mm$,  since
 sections of $\wedge^{*}A$ admit  nice operations such as the
exterior derivative and the Schouten bracket.

\begin{thm} ( \cite{Liu})
\label{th:pair}
 Let $(A,A^*)$ be a Lie bialgebroid, $L \subset A \oplus A^*$ a maximal
 isotropic subbundle defined by  a characteristic pair $\ch$
as in Equation (\ref{eq:L}).
 Then $L$ is a Dirac structure iff the following three conditions hold:
\begin{enumerate}
\item $D\subseteq A$ is a Lie subalgebroid.
\item $H$ satisfies the Maurer-Cartan type equation (mod D):
\begin{equation}
\label{MD}
d_{*} H+\half [H, H] \equiv 0,  (mod D).
\end{equation}
\item $ \gm (D^{\perp})$ is closed under the bracket
 $[\cdot, \cdot ] + [\cdot, \cdot ]_{H}$,
 where $[\cdot, \cdot ]_{H}$ is given by Equation (\ref{eq*}). I.e.,
\begin{equation}
\label{Dp}
 [\xi, \eta ] + [\xi, \eta ]_{H}  \, \in \Gamma (D^{\perp}), \ \ \ \
 \forall \xi , \eta \in \Gamma (D^{\perp}).
\end{equation}
\end{enumerate}
\end{thm}

Dirac structures are important in the construction of
Lie bialgebroids and Poisson homogeneous spaces.
For details, readers may consult the  references \cite{LWX} and \cite{LWX1}.

Finally,   note that  we may also work over $\complex$
when $M$ is a complex manifold. In this case, we just need
to replace smooth functions
by holomorphic functions, and smooth
sections by  holomorphic sections etc., and all the results
above will also hold.  In the
sequel, we will mainly work with complex Lie algebroids.
Even though one normally  works with sheaf of local sections 
when dealing with complex Lie algebroids since there may
not exist many  global sections. However, in the case below,
we can still avoid using sheaf  since we are working on 
an open subset $U$ of $\complex^n$.


\section{Twists of the standard $r$-matrix} 

Dynamical $r$-matrices have appeared  in various contexts \cite{Avan} \cite{EV}
\cite{F} \cite{Lu}.  In this section, we will show how a 
 dynamical $r$-matrix arises
 naturally as a twist of the standard classical $r$-matrix in the category of
Lie bialgebroids.

Let $\frakg$ be a simple Lie algebra over $\complex$   with a fixed
Cartan subalgbra $\f{h}$  and a root space decomposition:
\begin{equation}
\label{rsd}
\rsd =   {\frak n_+} \oplus \frakh  \oplus {\frak n_-},
\end{equation}
where $\f{n_{\pm}} = \sum_{\alpha \in \Delta_{\pm}}{\frak g_{\alpha}}$.
Let $< \cdot , \cdot >$ denote the Killing form on $\f{g}$ and
$E_{\alpha }\in \frakg_{\alpha }$ 
such that\\
$ <E_{\alpha}, E_{-\alpha}> = 1$.
Then  the standard classical r-matrix $r_0$ takes the form:
\begin{equation}
\label{r_0}
r_0 = \sum_{\alpha \in \Delta_{+}} E_{\alpha}\wedge E_{-\alpha}.
\end{equation}
Let $h_{\al} =  [E_{\alpha}, E_{-\alpha}] \in \f{h}$ for $\al \in
\Delta_{+}$  and  $h_{i} = h_{\al_{i}} $
for  simple roots  $\alpha_{i}$, $i=1, \cdots n$.  Then
 $\{h_1, \cdots ,  h_n\}$  forms a basis
 of $\f{h}$. Let  $\{h^{*}_1, \cdots , h^{*}_n \}$ be
its   dual basis,  which in turn induces a coordinate
system $(\lambda_{1}, \cdots , \lambda_{n})$ of
  $\frakh^*$, i.e.,
$\lambda = \sum \lambda_i h^{*}_i$ ,\,  $\forall \lambda \in \frakg^*$.

Now let   $U \subset {\frak h^*}$ be a connected open subset.
Consider the gauge  Lie algebroid:
\begin{equation}
\label{A1}
A = TU \times {\frak g} \cong U \times ({\frak h^*} \oplus {\frak g}).
\end{equation}
Set
\begin{equation}
\label{theta}
\theta = \sum_{i=1}^{n} h_{i} \wedge {\plam} .
\end{equation}


Clearly  $\theta $ can be considered as a constant
section of $  \wedge^{2} A$.
Equip $A^*\cong T^{*}U\times \frakg^*$ with the product
Lie algebroid, where $T^{*}U$ is the  trivial Lie algebroid
and $\frakg^*$ is the dual  Lie algebra induced by $r_0$.
Then
 $(A, A^*, r_0 )$ is a coboundary  Lie bialgebroid 
 according to  Proposition \ref{pro:gauge}.

\begin{thm}
\label{main1}
Let  $ \T  :U\lon  \wedge^{2}\f{g} $ be a holomorphic
functions considered as a section of $\wedge^{2}A$.
 Then   $\theta + \T$ is  a Hamiltonian
operator  of the Lie bialgebroid $(A, A^{*}, r_{0})$
iff $r = r_0 + \T$, the twist of $r_0$ by $\T$, is a dynamical
 $r$-matrix. 
\end{thm}
\pf $\theta + \T$ is a Hamiltonian
operator  iff it satisfies the
  Maurer-Cartan type  equation  (see Equation (\ref{eq:M})):
\begin{equation}
\label{mc1}
     d_*(\theta + \T) + \half [\theta + \T, \theta + \T] =0.
\end{equation}
By definition, $d_*(\theta + \T) = [ r_0, \theta + \T]$.
 Since $r_0$ is $\frakh$-invariant and  independent
of $\lambda$,
we have $[ r_0, \theta] = 0$.
It is also easy to see that $ [\theta, \theta ] = 0$,  and  
 $[\theta, \T] = \sum (h_i \wedge \plamm{\T}  + [ h_i, \T ] \wedge \plam)$.
Thus Equation (\ref{mc1}) becomes:
\begin{equation}
\label{mc2}
\begin{array}{lll}
-\sum [ h_i, \T ] \wedge \plam & =&  (\sum h_i \wedge \plamm{\T}) +
[r_0, \T]+ \half [\T,\T] \\
    & =&  (\sum h_i \wedge \plamm{(
    r_0 + \T)}) + \half [r_0 + \T, r_0 + \T] - \half [r_0, r_0] \\
    & =&  Alt(dr) + \half [r, r]  - \half [r_0, r_0].
\end{array}
\end{equation}
Now the left side of Equation (\ref{mc2}) belongs to
$\gm (\frakg \wedge  \frakg \wedge TU )$,  whereas
the right hand  side is a section of the subbundle
 $ \wedge^3 (U \times \f{g}) $. 
Thus both sides have to be zero identically.
This implies that 
$ [ h_i, \T ]=0$, $\forall i$, i.e.,
 $\T$ is $\f{h}$-invariant,  and  $r$ satisfies the modified
  CDYBE (\ref{rr})  since $\half [r_0, r_0] = [\Omega^{12},\Omega^{23}]$.
\qed

Now assume that $r=r_{0}+\T$ is  a dynamical r-matrix. Therefore
$\theta +\T$ is a Hamiltonian operator so that
its graph $\gm_{\theta+\T}$ is a Dirac structure of  the double
of $(A, A^{*}, r_{0})$.
Clearly, $\gm_{\theta+\T}$ is transversal to
$A$, so $(A,  \Gamma_{\theta +T})$ is
a Lie bialgebroid according to
 Theorem 2.6 in \cite{LWX}.
 In fact, it is simple to see that the Lie algebroid $\gm_{\theta +\T}$ is
isomorphic to $A^*$ with a twisted bracket defined by
the new  $r$-matrix  $\Lambda :=
\theta + \T + r_0 = \theta +r$, so
$(A, A^{*}, \Lambda )$ is also a coboundary
Lie bialgebroid. Thus, we have proved the following
result of  Bangoura and Kosmann-Schwarzbach \cite{BK-S}:
 
\begin{cor}
\cite{BK-S}
Let  $r(\lambda ):U\lon \wedge^{2} \frakg $ be a holomorphic
function.  Then $\Lambda =\theta +r(\lambda ) \in \gm (\wedge^{2}A )$
 defines a  coboundary Lie bialgebroid  iff $r(\lambda )$
 is a dynamical r-matrix.
\end{cor}

It is not difficult to see that this  Lie bialgebroid
is the Lie  bialgebroid corresponding to the
dynamical Poisson groupoid constructed by Etingof and Varchenko \cite{EV}.
The following conclusion follows immediately from the construction.
\begin{thm}
Let  $r(\lambda ) $ be a dynamical
$r$-matrix,  and  $\Lambda =\theta +r(\lambda ) $ the twisted $r$-matrix.
Then, as a Courant algebroid,
  the double of the coboundary Lie bialgebroid
$(A, A^{*}, \Lambda )$  
is  isomorphic to the double of the untwisted Lie bialgebroid
$(A, A^{*}, r_0 )$. 
\end{thm}

It is simple to see  that a function
 $\T : U  \lon  \wedge^{2} \f{g}$
  is $\frakh$-invariant  iff it can be  splitted into 
 two terms: $\T = \omega + \T_0$ , where 
\begin{equation}
\label{T} 
 \omega = \sum_{ij}  \omega ^{ij} (\lambda )h_{i} \wedge h_{j}, \ \ \ \
 \mbox{and } \ \ \ \ \T_0 = \sum_{\alpha \in \Delta_{+}}
  \T_{\alpha} (\lambda )E_{\alpha}\wedge E_{-\alpha}.
\end{equation}

\begin{pro}
Let  $\T$  be given as above. Then $ \theta + \T$ is  a Hamiltonian
operator iff
\begin{enumerate}
\item $ \T_0$ is a Hamiltonian operator; and
\item    $\omega $ is a  closed 2-form on $ U$.
\end{enumerate}
\end{pro}
\pf The Maurer-Cartan equation for  $\theta + \T_0 + \omega$ takes the form:

\be
0 &=&d_*(\theta + \T_0 + \omega) +\half [\theta + \T_0 + \omega,
\theta + \T_0 + \omega] \\ 
 &=& d_*(\theta + \T_0 ) +\half [\theta + \T_0 , \theta + \T_0 ]
 + [\theta , \omega]. 
\ee

Note that, on the right hand  side of the  equation, the only  term
in $ \wedge^3 \f{h}$ is 
$$[\theta , \omega] = \sum h_i \wedge \plamm{\omega} = d\omega .$$
So  the equation holds iff 
$$
d_*(\theta + \T_0 ) +\half [\theta + \T_0 , \theta + \T_0 ] = 0,
\ \ \ \ \  and \ \ \ \ d\omega =0.
$$
Thus the proposition is proved.
\qed

In the terminology of Etingof and Varchenko, $\T$ and $\T_0$ are called
gauge equivalent, and $\omega$ is a gauge term. In fact, for
most purposes we may assume that $\omega =0$.

Finally, note that for any fixed $\lambda \in U$,
$\T (\lambda )\in \wedge^{2}\frakg$ is generally not a Hamiltonian
operator for the Lie bialgebra $(\frakg , \frakg^* , r_{0})$.
In fact, it is easy to see that  $r=r_{0}+\T$ is a dynamical $r$-
matrix iff 
\begin{equation}
\label{eq:t}
[r_{0}, \T ] + \half [\T, \T ] + Alt (d\T) = 0.
\end{equation}
Thus,
\be
 d_{*}\T(\lambda )  +\half [\T(\lambda ), \T(\lambda )] 
  &=&[r_{0}, \T(\lam)]+\half [\T(\lam), \T(\lam) ] \\
  &=&([r_{0}, \T ]+\half [\T , \T ])(\lam ) \\
  &=& -Alt (d\T) (\lam). 
\ee
So  $\T (\lambda )$ is a Hamiltonian operator  iff $\lambda$
is a critical point of $\T$
(we will see in Section 4  that this
is equivalent to $ \T \equiv 0$ on $U$).
Hence  $-Alt(d\T) (\lam )$ measures the failure of the graph of
 $\T(\lam)^{\#}: \f{g^*} \lon \f{g} $
 being a Lagrangian subalgebra.
In terms of  Drinfel'd \cite{dr:quasi}, $\T (\lambda )$  is
a family of twists, which defines a
family of quasi-Lie bialgebras $(\f{g}, \delta (\lambda) , \phi (\lambda ))$.
Here $\delta (\lambda ):\frakg \lon \wedge^{2}\frakg$ is 
given by $\delta (\lambda )(x)=[r_{0}+\T(\lambda ) , x]$ and
$\phi (\lambda )=-Alt(d\T) (\lam )\in \wedge^{3}\frakg$.
This family of quasi Lie bialgebras
is  the classical limit  of  the quasi-Hopf algebras
studied by Fronsdal \cite{Fron}, Arnaudon et. al.  \cite{Arnaudon} and 
Jimbo et. al. \cite{Jimbo} connected 
with  quantum dynamical $R$-matrices (see also \cite{Xu}).

\section{Construction of Dirac structures}

In the previous  section, we have already established
a simple connection between dynamical $r$-matrices
and Dirac structures. The purpose of this section is
to give an explicit construction of those Dirac
structures.

As in Section  3,  assume that $\frakg$ is a simple
Lie algebra with Killing form $<\cdot,\cdot>$, and $r_{0}=
\sum_{\alpha \in \Delta_{+}}E_{\alpha}\wedge E_{-\alpha}$
is  the standard $r$-matrix. By identifying $\frakg^*$ with
$\frakg$ using the Killing form, the bracket on $\frakg^*$
is given by:
 $$ [X,Y]_{R} = [R X,Y]+[X, R Y], \ \ \ \ \forall X, Y \in {\frak g},$$
where
$ R= \pi_{+} - \pi_{-}$, and
  $\pi_{\pm} : \f{g} \lon \f{n_{\pm}}$ are the
natural  projections with respect to
 the Gauss decomposition
$\f{g} =\f{n_{+}} \oplus \f{h} \oplus \f{n_{-}}$ as  in Equation (\ref{rsd}).
It is well-known that the double of the Lie bialgebra $(\frakg, \frakg^{*})$
can be  identified with   the direct sum Lie algebra
${\frak d} = {\frak g} \oplus {\frak g}$,  while the corresponding
 invariant non-degenerate  bilinear form is:
\[
( (X_1, Y_1), \, (X_2, Y_2) ) \, = \, {1 \over 2}
(< Y_1 , Y_2 > \, - \, < X_1, X_2 >), \ \ \forall X_{1}, X_{2}, Y_{1}, Y_{2}
\in \frakg.
\] 
Here ${\frak g}$ is identified with the diagonal,
while  ${\frak g^*}$ is identified with the subalgebra:
\[
\{(X_{-}+ h, X_{+} -h)| ~  \forall X_{\pm} \in {\frak n_{\pm}},h \in {\frak h} \}.
\]
Thus the corresponding Courant algebroid,  as the
double of the Lie bialgebroid $(A, A^*, r_0)$, is a trivial
vector bundle, which can be expressed as:
$$E= A \oplus A^*  \cong (TU \oplus T^*U) \times  \f{d}
\cong U \times ( \f{h} \oplus \f{h} \oplus \f{g} \oplus\f{g}). $$
Consequently, a section of $E$ can be considered as a vector-valued
function on $U$ with value in 
 $ \f{h} \oplus \f{h} \oplus \f{g} \oplus\f{g}$, which is denoted by
$\4{\xi (\lambda )}{\eta (\lambda )}{X (\lambda )}{Y (\lambda )} $. Here
$\xi (\lambda ), \eta (\lambda ) $ are  $\f{h}$-valued
functions on $U$, 
and $X(\lambda ), Y (\lambda )$ are $  \f{g} $-valued
functions on $U$.
The inner-product $(\cdot,\cdot)$ on $E$ is given by
\begin{equation}
\label{product}
\begin{array}{lll}
 (( \xi_1, \eta_1 ; X_1, Y_1), ( \xi_2, \eta_2 ; X_2, Y_2))& = &
{1 \over 2}(< \xi_1, \eta_2 > \, + \, < \eta_1, \xi_2 >) \\
&& +{1 \over 4}  (< Y_1, Y_2 > \, - \, < X_1, X_2 >).
\end{array}
\end{equation}
 Then as subbundles  of $E$, $A$ and $A^*$ are  given by
\begin{eqnarray}
\label{A}
A & \cong& U\times \{\4 k{0}{X}{X} ~|~ \,\forall  k \in \f{h}, X \in \f{g}\},\ \ \ \mbox{ and}\\
A^* & \cong & U\times \{\4 0{h}{X_{-}+k}{X_{+}-k} \,
 ~|~ \, \forall  h, k \in \f{h}, X_{\pm} \in \f{n_{\pm}}\}.
\label{A*}
\end{eqnarray}
As for the bracket of $\gm (E)$, it admits   a simple form
for constant sections:
\begin{equation}
\label{form0}
[ \4 {\xi_{1}}{\eta_{1}}{X_1}{Y_1}, \,  \4 {\xi_{2}}{\eta_{2}}{X_2}{Y_2}]
\, = \, \4 0{0}{[X_1,X_2]}{[Y_1 ,Y_2]}.  
\end{equation}
For general sections, the formula is much involved.
The following are two special cases corresponding  to Equations
(\ref{eq:fg}) and (\ref{eq:fe}), which are needed
in the future:
\begin{equation}
\label{form}
[ (0, 0; 0,fX), \4 0{0}{0}{gY}]=\4 0{{1 \over 4}(gdf-fdg)<X,Y>}{0}{fg[X,Y]},
\end{equation}
and
\begin{equation}
\label{form1}
[ (h^{*}_i, h_j ; 0,0), \4 0{0}{fX}
{gY}] \,= \, \4 0{0}{{\partial f \over \partial \lambda_i}X}
{{\partial g \over \partial \lambda_i}Y},
 \ \ \ \ \forall f, g \in \cinf (U), \, X,Y \in {\frak g}.
\end{equation}

Next we need to describe the graph of  $\theta^{\#} +\T^{\#}: A^* \lon A$.
For simplicity  we assume that $\T$ is given by Equation
 (\ref{T}) with $\omega =0$.  
Set 
$$ d = \{ \4 {k}{h}{h+k}{h-k} \, ~|~
 \, \forall h, k \in \f{h} \,\} \subset
 \f{h} \oplus \f{h} \oplus \f{g} \oplus \f{g}.$$
And for each $\lambda\in U$, define
\begin{equation}
\label{B}
B(\lambda) =\{ \4 0{0}{(r^{\#}(\lambda)-1)X}{(r^{\#}(\lambda)+ 1)X}
 \,|\, \forall X \in {\frak n_{\pm}}\, \}, 
\end{equation} 
where $r(\lambda )=r_{0}+\T(\lambda )$ as in  Theorem \ref{main1}.

\begin{lem}
\label{lem:graph}
As a subbundle of $E$,  the graph   of  $\theta^{\#} +\T^{\#}:
 A^* \lon A$ is $L=\cup_{\lambda \in U}L(\lambda )$, where
$$ L(\lambda )  =d  \oplus B (\lambda ).$$ 
\end{lem}
\pf Using the identification as in Equations (\ref{A}) and (\ref{A*}),
we need to
compute the image $(\theta^{\#}+\T^{\#}) \4 0{h}{X_- + k}{X_+ -k}$
at each $\lambda \in U$. Now
\be
\theta^{\#} \4 0{h}{X_{-} + k}{X_{+} -k} &=& \4 {k}{0}{h}{h}, \
\ \mbox{ and }\\
\T^{\#}  \4 0{h}{X_{-} + k}{X_{+} -k}&=&
 \half \4 0{0}{\T^{\#} X_+}{\T^{\#} X_+}  -  \half \4 0{0}{\T^{\#} X_-}
 {\T^{\#} X_-}.
\ee

Therefore,
\be
&&(\theta^{\#}+\T^{\#}) \4 0{h}{X_- + k}{X_+ -k}
+ \4 0{h}{X_- + k}{X_+ -k} \\
&=& \theta^{\#} \4 0{h}{k}{-k} + \4 0{h}{k}{-k}\\
&&+ \T^{\#} \4 0{0}{X_- }{X_+ } + \4 0{0}{X_-}{X_+} .
\ee

It is easy to see that
$$
 \theta^{\#} \4 0{h}{k}{-k} + \4 0{h}{k}{-k}=
\4 k{h}{h+k}{h-k}  \in d.
$$
And
\be
&&\T^{\#} \4 0{0}{X_- }{X_+ } + \4 0{0}{X_-}{X_+}\\
&=& \half \4 0{0}{\T^{\#}X_+}{(\T^{\#}+2)X_+} -
 \half \4 0{0}{(\T^{\#}-2)X_-}{\T^{\#}X_-}\\
&=&\half \4 0{0}{(r^{\#}-1)X_+}{(r^{\#}+1)X_+}
 - \half \4 0{0}{(r^{\#}-1)X_-}{(r^{\#}+1)X_-} \in B(\lambda ),
\ee
where we have used the  fact that $r^{\#}| _{\f{n}{\pm}} = \T^{\#}  \pm 1$. 
This concludes the proof of the lemma.  \qed

 For any  $\lambda \in U$, consider the decomposition:
\begin{equation}
\label{kn1}
\f{n_{\pm}}= \f{k_{\pm}}(\lambda) \oplus \f{n^{\circ}_{\pm}}(\lambda),
\end{equation}
where                     
\begin{equation}
\label{ker}
\f{k_{\pm}} (\lambda ) = ker \T^{\#} (\lambda) \cap \f{n_{\pm}}
= span_{\complex}
 \{ E_{\pm \alpha} \, | \T_{\alpha}(\lambda) =0, \alpha \in \Delta_{+} \},
\end{equation}
and
\begin{equation}
\f{n^{\circ}_{\pm}} (\lambda )=  span_{\complex} \{ E_{\pm \alpha} \, |
\, \T_{\alpha}(\lambda) \neq 0, \alpha \in \Delta_{+} \}.
\end{equation}
 Then we can  rewrite $B(\lambda)$ as follows:
\begin{equation}
\label{B1}
B(\lambda)= span_{\complex}\{\4 0{0}{X}{\varphi(\lambda) X },
 \4 0{0}{Y_{-}}{Y_{+}} \,|\, \forall X
 \in {\frak n}^{\circ}_{\pm}(\lambda ),  Y_{\pm}
  \in{\frakk_{\pm}}({\lambda})\},
\end{equation}
where
$$\varphi(\lambda) = { {r^{\#}(\lambda)+1} \over {r^{\#}(\lambda)-1}}
\,: \ \ \ \ \f{n^{\circ}_{\pm}} (\lambda ) \lon  \f{n^{\circ}_{\pm}} (\lambda )
$$  is the  Cayley transformation of  the linear operator
$r^{\#}(\lambda) \,|_{ \f{n}^{\circ}_{\pm}(\lambda)}$.
Consequently, $L$ can be written as:
\begin{equation}
\label{eq:Lr}
\begin{array}{lll}
L(\lam )&=&
span_{\complex} \{ \4 {k}{0}{k}{-k}, \4 0{h}{h}{h}, 
 \,\4 0{0}{X}{\varphi  (\lambda )X },  \4 0{0}{Y_-}{Y_+}~  
|\\
&& \mbox{} \hspace{0.6in} \forall  h, k \in \f{h}, 
X \in \f{n^{\circ}_{\pm}}(\lambda),  Y_{\pm} \in
k_{\pm}(\lambda) \}.
\end{array}
\end{equation}

\begin{lem}
\label{le:S}
Assume that    $L\subset E $  is  a  Dirac structure, then 
\begin{enumerate}
\item both $ \f{k_{\pm}}(\lambda)$ and  $\f{n^{\circ}_{\pm}}(\lambda) $
  are independent of $ \lambda \in U$ (for
simplicity, we  denote  them by $ \f{k_{\pm}}$ and  $\f{n^{\circ}_{\pm}}$
respectively);
\item $\f{n^{\circ}_{\pm}}$ are subalgebras of $\f{n_{\pm}}$;
\item  $\f{k_{\pm}}$ are ideals of   $\f{n_{\pm}}$.
\end{enumerate}  
\end{lem}
\pf According to  Theorem \ref{main1},
  $r_{0} + \T$ is  a dynamical $r$-matrix. By Equation 
(\ref{eq:t}), we have
\be
0 &=& [r_0 , \T ] +\half [\T ,  \T ]  + Alt (d \T )  \\
 &=& \sum_{i}
  \plamm {\T}  \wedge  h_i +
\sum_{\alpha , \beta \in \Delta_{+}} [(\half \T_{\alpha}+ 1)
E_{\alpha}\wedge E_{-\alpha}, \ \  \T_{\beta}E_{\beta}\wedge E_{-\beta}] \\
&=& 
\sum \limits_{\alpha \in \Delta_{+}} 
\sum \limits_ {i} \plamm {\T_{\alpha}} 
E_{\alpha}\wedge E_{-\alpha} \wedge h_{i}
~ + ~\sum \limits_{\alpha , \beta \in \Delta_{+}} (\half \T_{\alpha}+ 1) \T_{\beta}
[E_{\alpha}\wedge E_{-\alpha},  \  E_{\beta}\wedge E_{-\beta}].
\ee 
Since $ [E_{\alpha},\, E_{-\alpha}] = h_{\al}
= \sum <\al, h^*_i>h_i$
 for  any $ \alpha \in \Delta_{+}$,
the coefficient of the term $ E_{\alpha}\wedge E_{-\alpha} \wedge h_{i}$
in the above equation is
$\plamm{\T_{\alpha} } - <\al, h^*_i>(\T_{\alpha}+ 2)\T_{\al}$.
This implies that $\T_{\alpha}$ satisfies the
following system of first-order differential equations: 
$$  \plamm { \T_{\alpha}} - <\al, h^*_i> (\T_{\alpha}+ 2) \T_{\alpha}
 = 0 , \ \ \ \  \forall \alpha \in \Delta_{+},
\, (i = 1, \cdot \cdot \cdot, n).
$$
Thus  if $\T_{\alpha}(\lambda_{0}) = 0$ for some
$\lambda_0 \in U$, then $\T_{\alpha} \equiv 0$ on $U$.
This   is equivalent to  that 
$\f{k_{\pm}} (\lambda ) = ker \T(\lambda) \cap \f{n_{\pm}}$
are independent of $ \lambda \in U$. This proves the first statement.

For the second statement, note that since $r(\lambda )$ 
is $\frakh$-invariant,  $\va  (\lambda )$ commutes with $ad_{\frakh}$.
 Thus   $\va E_{\al} = \va_{\al}E_{\al}$ for some function $\va_{\al}:
U\lon \complex$, $\forall \alpha \in \f{n^{\circ}}_+$.
For any  $\alpha, \beta \in \f{n^{\circ}}_+$,  since
$\4 0{0}{E_{\al}}{\va_{\al}E_{\al}}
,\4 0{0}{E_{\beta}}{\va_{\beta}E_{\beta}} \in \Gamma (L)$,
their commutator  belongs to $\gm (L)$ as well.

On the other hand, it is clear that
\be
[\4 0{0}{E_{\al}}{\va_{\al}E_{\al}}, \4 0{0}{E_{\beta}}{\va_{\beta}E_{\beta}}]
 & = &\4 0{0}{[E_{\al}, E_{\beta}]}{\va_{\al}\va_{\beta} [E_{\al}, E_{\beta} ]} \\
 && +  \4 0{\frac{1}{4}(\va_{\beta}d\va_{\al} - \va_{\al} d \va_{\beta} )
 <E_{\al}, E_{\beta} > }{0}{0} \\
 & =& N_{\al,\beta}  \4 0{0}{E_{\al+ \beta}}{ \va_{\al}\va_
 {\beta}E_{\al+ \beta}}.
\ee
Here, in the last equality, 
 we used the fact that $ <E_{\al}, E_{\beta} > = 0$ whenever
$\al \neq \beta$.
According to  Equation (\ref{B1}), we conclude  that
 $E_{\alpha + \beta} \in \f{n^{\circ}}_+$ whenever
  $N_{\al ,\beta} \neq 0$, i.e, $\alpha + \beta \in \Delta_{+}$.
 This means  that $\f{n^{\circ}}_{+}$ is  a  Lie  subalgebra of $\f{n_{+}}$
and
\begin{equation}
\label{a+b}
  \va_{\al} \va_{\beta} \, = \, \va_{\al+ \beta }, \ \ \ \ \forall 
E_{\alpha}, E_{\beta} \in \f{n^{\circ}}_+  \ \ \ \ \mbox{such that} \ \ \ \
 \alpha + \beta \in \Delta_{+}.
\end{equation}
Similarly we can prove  that $\f{n^{\circ}}_{-}$ is  a
Lie  subalgebra of $\f{n_{-}}$.

For the third statement, let $X_+, Y_{+} \in \f{k_{+}}$,
and $E_{\al} \in  \f{n^{\circ}}_+$.
As constant sections of $\gm (L)$, we have
$$
[\4 0{0}{0}{X_{+}}, \, \4 0{0}{0}{Y_{+}}] \, = \,\4 0{0}{0}{[X_+, Y_{+}]}
\, \in \gm (L),
$$
which means that $ [X_+, Y_{+}] \in  \f{k_{+}}$. Moreover,
$$
[ \4 0{0}{E_{\al}}{\va_{\al} E_{\al}} , \4 0{0}{0}{Y_+}]
 = \4 0{0}{0}{\va_{\al} [E_{\al}, Y_+] }  \in \gm (L).
$$
This implies  that $[E_{\al}, Y_+] \in \f{k_{+}}$. Thus 
 $\f{k_{+}}$ is an ideal of $\f{n_{+}}$ since
$\f{n_{+}} = \f{k_{+}} \oplus   \f{n^{\circ}}_+ $.  
 Similarly, $\f{k_{-}}$ is an ideal of $\f{n_{-}}$.  \qed

Below we will see that  any decomposition
$\f{n_{\pm}}= \f{k_{\pm}} \oplus \f{n^{\circ}_{\pm}}$ satisfying 
Properties (2)-(3) in Lemma \ref{le:S} corresponds
to a  subset $S$ of simple roots. More precisely,
given a decomposition 
$\f{n_{\pm}}= \f{k_{\pm}} \oplus \f{n^{\circ}_{\pm}}$,
let $S$ be  the subset of those simple roots $ \al_i $ such that
$E_{\al_i} \in \f{n^{\circ}_+}$.
Define a  subset of positive roots as follows:
\begin{equation}
\label{S1}
[S] = \{ \al \in \Delta_+ ~|~  \al = \sum_{ \al_i \in S}n_i\al_i, \,\,
n_i \geq 0 \, \}.
\end{equation}
Since any positive (negative) root can be expressed as positive (negative)
linear combination of simple roots, we have
\begin{pro}
\label{kn}
Assume that   $\f{n_{\pm}}= \f{k_{\pm}} \oplus \f{n^{\circ}_{\pm}}$ 
is a  decomposition satisfying  Properties (2)-(3) in Lemma \ref{le:S}.
Then,
\begin{equation}
\label{eq:n}
 \f{n^{\circ}}_{\pm}  =  span_{\complex} \{E_{\pm \al}, \,  \al \in  [S] \},
\end{equation}
i.e., $\{ E_{\pm \al_i}| \al_{i} \in S \}$  are Lie algebraic generators of
$\f{n^{\circ}}_{\pm}$. Consequently, 
\begin{equation}
\label{eq:k}
\f{k}_{\pm}  =  span_{\complex}
 \{E_{\pm \al}, \, \al \in \Delta_{+} \backslash  [S] \}.
\end{equation}
Conversely, given any subset $S$ of simple roots, the
corresponding $ \f{k_{\pm}}$ and $ \f{n^{\circ}}_{\pm} $ defined  by
Equations (\ref{eq:n}) and (\ref{eq:k})
above satisfy  Properties (2)-(3) in  Lemma \ref{le:S}.
\end{pro}

Now we are ready to  prove the main theorem of this section.

\begin{thm}
\label{construction}
Let $S$  be a subset of simple roots
 with corresponding $ \f{k_{\pm}}$ and $ \f{n^{\circ}}_{\pm} $ 
defined as in Proposition  \ref{kn},
and $L$  a  subbundle of $E$ defined by Equation (\ref{eq:Lr}),
 where $\va (\lambda ), \forall
\lambda \in U$,  is a linear operator on $n^{\circ}_{\pm}$.
Then $L$ is a Dirac structure iff there exists
 some  $\lambda_{0} \in \f{h}$ such that
  $ \varphi(\lambda)  =
\,  Ad_{ e^{2( \lambda +\lambda_{0})}}$.
\end{thm}
\pf We shall divide the proof into four steps.

{\it Step 1}.  It follows from Equations (\ref{form0}),
(\ref{form}) and (\ref{form1})
that for any    $ h \in \f{h}$ and
\,$X \in \f{n^{\circ}}_{\pm}$, 
$$
[ \4 {0}{h}{h}{h}, \, \4 0{0}{X}{\varphi X} \,] =
 \4 0{0}{[ h,X]}{[h, \varphi X]}.
$$
It  is still in $L$ iff
$$ [h, \varphi X] = \va [h, X],$$
which is equivalent  to that $\varphi$ commutes with 
$ad_{{\frak h}}$. Therefore,
 $\varphi E_{\alpha} =\varphi_{\alpha}E_{\alpha}$,
  $\forall \alpha \in \pm [S]$,  i.e., $ E_{\al} \in n^{\circ}_{\pm}$, 
where $ \varphi_{\alpha}$ is a complex valued function on $  U$.

{\it Step 2}.  Suppose  that $\varphi$ commutes with  $ad_{{\frak h}}$.
Then
$\forall  i = 1, \cdot \cdot \cdot, n$ and $ E_{\al} \in \f{n^{\circ}}_{\pm}$,
both $\4 {h^{*}_i}{0}{h^{*}_i}{-h^{*}_i}$ and
$\4 0{0}{E_{\al}}{\varphi_{\al}E_{\al}}$ are sections of  $L$.
By Equations (\ref{form0}),  (\ref{form}) and (\ref{form1})

\be
[\4 {h^{*}_i}{0}{h^{*}_i}{-h^{*}_i},\4 0{0}{E_{\al}}{\varphi_{\al}E_{\al}}]
&=&\4 0{0}{0}{\plamm{\va_{\al}} E_{\al}}+
\4 0{0}{[h^{*}_i, E_{\al}]}{- \va_{\al}[h^{*}_i, E_{\al}] } \\
&=& \4 0{0}{<\al, h^{*}_i>E_{\al}}{(\plamm{\va_{\al}}
 - <\al, h^{*}_i>\va_{\al})E_{\al} }.
\ee
It is still in $\gm (L)$  iff 
$$\plamm{\va_{\al}} = 2 <\al, h^{*}_i>\va_{\al}  \ \ \ \ \Longleftrightarrow \ \ \ \
  \varphi_{\alpha}(\lambda) = C_{\alpha} e^{2<\alpha, \lambda>},$$
where $C_{\al}$ are certain constants and $\lambda = \sum \lambda_i h^{*}_i$.

{\it Step 3}.  Suppose that $\varphi_{\alpha}(\lambda) = C_{\alpha}
e^{2<\alpha, \lambda>}$. Next we show that  $C_{\alpha}$
satisfy the following relations:
\begin{equation}
\label{constant}
 C_{-\alpha}=  C_{\alpha}^{-1},  \ \ \ \ C_{\alpha+\beta}
   = C_{\alpha} C_{\beta}, \ \ \ \ 
\forall \alpha, \beta \in \pm [S], \ \ \mbox{ whenever }
 \alpha +\beta \in \Delta.
\end{equation}
When $\beta \neq - \al$ the conclusion follows  from
Equation (\ref{a+b}), where only  a special case:
both $\al$ and $\beta$ being positive roots, is discussed. However,
 the general situation can also be
 easily   checked using the  fact that $< E_{\al}, E_{\beta}>=0$.
Now assume that  $\beta = -\al$. Then, by Equation (\ref{form}), we have 
\be
&&[ (0, 0; E_{\alpha},\varphi_{\alpha} E_{\alpha}),
 \4 0{0}{E_{-\alpha}}{\varphi_{-\alpha}E_{-\alpha}}] \\
&=&[ (0, 0; E_{\alpha}, \, e^{2<\alpha, \lambda>} C_{\alpha} E_{\alpha}), \,
 \4 0{0}{E_{-\alpha}}{e^{-2<\alpha, \lambda>}C_{-\alpha}E_{-\alpha}}] \\
&=&\4 0{ {1 \over 4 }(
e^{-2<\alpha, \lambda>} de^{2<\alpha, \lambda>} -
e^{2<\alpha, \lambda>} de^{-2<\alpha, \lambda>} )
<E_{\alpha},E_{-\alpha}>}
{[E_{\alpha},E_{-\alpha}]}
{C_{\al}C_{-\al}[E_{\alpha},E_{-\alpha}]} \\
&=&\4 0{\, \sum_{i}
{{\partial <\alpha, \lambda>} \over {\partial \lambda_i}} d\lam_{i} \,}
{[E_{\alpha},E_{-\alpha}]}
{C_{\al}C_{-\al}[E_{\alpha},E_{-\alpha}]} \\
&=& \4 0{\sum <\alpha, h^{*}_i>h_i}{h_{\alpha}}{C_{\al}C_{-\al} h_{\alpha}} \\
&=& \4 0{h_{\alpha}}{h_{\alpha}}{C_{\al}C_{-\al} h_{\alpha}},
\ee
where  we used the facts that
\begin{equation}
\label{facts}
 d\lam_{i} = h_i, \ \   \lam = \sum  \lam_i h^{*}_i, \ \
h_{\al} = \sum <\al, h^{*}_i>h_i \ \  and  \ \ < E_{\alpha},E_{-\alpha}>=1.
\end{equation}
Obviously the commutator is still in $\gm (L)$  iff $C_{\al}C_{-\al} =1$.
Thus, Equation (\ref{constant}) is proved.

Finally, it is not difficult to see that
Equation (\ref{constant}) implies  that there exists some
 $\lambda_{0} \in {\frak h}$ such that
$ C_{\alpha} = \,  e^{2<\alpha, \lambda_{0}>}$.
In fact,  we can take
$\lambda_{0} =  \half {\sum_{ \al_i \in S}  (ln \, C_{\al_i}) h^{*}_i} $.
Consequently, we have
\begin{equation}
\label{e}
 \varphi_{\al}(\lambda)  =
\, e^{2 <\al ,  \lambda +\lambda_{0} >}  \ \ \ \ \Longleftrightarrow
 \ \ \ \
\varphi(\lambda)  = \,  Ad_{ e^{2( \lambda +\lambda_{0})}}.
\end{equation}

Conversely,
if $\va (\lambda )= Ad_{ e^{2( \lambda +\lambda_{0})}}$,
  $L$ is maximal isotropic since  $ \varphi$ 
preserves the  Killing form $<\cdot, \cdot>$. Moreover,
$\gm (L)$ is closed, so $L$ is indeed
a Dirac structure. This concludes the proof.  \qed

\begin{cor}(\cite{EV})
\label{cor:ef}
A meomorphic function  $r: U \lon  \wedge^{2} \f{g}$ 
is a  dynamical $r$-matrix  iff
$r$ is of  the form:
\begin{equation}
\label{ef}
r(\lam ) = \omega + \,
\sum _{\alpha \in [S]} \,
\coth< \al, \lam + \lam_0 >
E_{\al} \wedge  E_{-\al } \,
+ \,
\sum _{\alpha \in \Delta_{+}\backslash [S] } \,
 E_{\al} \wedge  E_{-\al},
\end{equation}
where $\omega$ is a closed 2-form on $U$,  and $[S]$ is defined  by Equation
(\ref{S1}) for a subset $S$ of the simple roots.
\end{cor}
\pf
Let  $\T = r -r_0$. Then $\T$ is of the form:
\begin{equation}
\T = \omega + \, \sum _{\alpha \in \Delta_{+}} \,
\T_{\al}  E_{\al} \wedge  E_{-\al},
\end{equation}
where $\omega$ is a closed two-form on $U$.

According to Theorem \ref{main1},    
$r$ is a dynamical $r$-matrix iff
$\Gamma_{ \theta + \T } \subset  A \oplus A^* $ is a Dirac structure
of the  Lie bialgebroid $( A, A^* , r_0 )$. Without loss of generality,
assume that $\omega =0$.
According to  Theorem \ref{construction}, the latter amounts to
that there exists a subset  of simple roots $S$ with corresponding
$\f{n^{\circ}_{\pm}} $ and $\f{k_{\pm}}$ such that 
the Cayley transformation of $r^{\#}|_{\f{n^{\circ}}_{\pm}}$:
$\varphi(\lambda) = { {r^{\#}_1(\lambda)+1} \over {r^{\#}_1(\lambda)-1}}$
 has expression (\ref{e}),
for some fixed $\lam_{0} \in \f{h}$. This immediately  implies that
$$r(\lam )|_{\f{n^{\circ}}_{\pm}} = \sum _{\alpha \in [S]} \,
coth< \al, \lam + \lam_0 >
E_{\al} \wedge  E_{-\al},\ \ \ \ \mbox{ and}  \ \ \ \
r(\lam )|_{\f{k_{\pm}}} = 
\sum _{\alpha \in \Delta_{+}\backslash [S] } \,
 E_{\al} \wedge  E_{-\al}. $$
The conclusion thus follows.
\qed

\section{Lagrangian subalgebras and dynamical $r$-matrices}

In \cite{Ka},  Karolinsky classified all Lagrangian subalgebras $W_0$  of
the double of  the Lie bialgebra $( \f{g}, \f{g^*}, r_0)$
(by abuse of notation, in the sequel, 
 we will simply  say Lagrangian subalgebras of 
the Lie bialgebra $( \f{g}, \f{g^*}, r_0)$) in   terms of 
 the triples $(\f{u}^{-}, \, \f{u}^{+}, \, \va )$, where
$ \f{u}^{\pm}$ are two parabolic subalgebras of $\f{g}$ ,
$ \f{m} = \f{u}^{+} \cap \f{u}^{-}$ is a Levi subalgebra and
$\va$ is an inner automorphism of $\f{m}$. The following theorem shows
that such a classification can be reduced to a simpler form in  the
special case that $W_0 \cap \frakg =\frakh $.

\begin{pro}
There is a one-one correspondence  between
Lagrangian subalgebras  $ W_0 \subset \f{g} \oplus \f{g}$  with
 $W_0 \cap \f{g} = \f{h}$ and 
 pairs $( S, \lam_0 )$,  where $S$ is
 a subset  of simple roots and $\lam_0 \in \f{h^*}$.
\end{pro}
\pf Given  such a pair $( S, \lam_0 )$, define ${\frak n^{\circ}_{\pm}}$ and 
${\frak k_{\pm}}$ as in  Proposition \ref{kn} by:
\begin{equation}
\label{eq:n1}
 \f{n^{\circ}}_{\pm}  =  span_{\complex} \{E_{\pm \al}, \,  \al \in  [S] \},
\end{equation}
\begin{equation}
\label{eq:k1}
\f{k}_{\pm}  =  span_{\complex}\{E_{\pm \al}, \, \al \in \Delta_{+} \backslash  [S] \},
\end{equation}
where 
\begin{equation}
\label{S2}
[S] = \{ \al \in \Delta_+ ~|~  \al = \sum_{ \al_i \in S}n_i\al_i, \,\,
n_i \geq 0 \, \}.
\end{equation}
Let $W_{0}\subset \f{g} \oplus \f{g}$ be the subspace:
\begin{equation}
\label{eq:w0}
 W_{0}   = span_{\complex} \{ (h, h),
 \, (X , Ad_{e^{\lam_0}}X ), \, (Y_- , Y_+) \, ~
|\, \forall h \in \f{h},  \,X \in {\frak n^{\circ}_{\pm}}, \,
Y_{\pm} \in{\frak k_{\pm}} \}.
\end{equation}
One can check directly that $W_0$   is a Lagrangian subalgebra
and  $W_0 \cap \f{g} = \f{h}$. Here, as before,   the double
 is identified  with $\f{d} =\f{g} \oplus \f{g}$, 
whereas $ \f{g}$ is identified 
 with  the diagonal of $\f{d}$.

Conversely, as we know in Section 2, 
any Lagrangian subalgebra  of the double of a Lie bialgebra
arises  from a characteristic pair.
More precisely, given  a Lagrangian subalgebra  $ W_0 \subset \f{g} \oplus
\f{g}$ such that $W_0 \cap \f{g} = \f{h}$, 
there exists some $ J \in \f{g} \wedge \f{g}$ such that
\begin{equation}
\label{eq:J}
  W_0 = \{ \,X + J^{\#}\xi + \, \xi \,| \,  
 X \in \f{h} , \xi \in \f{h}^{\perp} \} = \f{h} \oplus \, graph(
J^{\#}|_{\f{h}^{\perp}}),
\end{equation}
i.e., $( \f{h}, J )$ is a characteristic pair of $W_0$.
However $J$ is not unique.
What we need here   is to choose an  $\frakh$-invariant
$J$.  For this purpose,
we notice that 
 $\forall X \in \f{h}$ and $\xi \in \f{h}^{\perp}$,
  \be
[ X,  J^{\#} \xi + \xi ] & = & [X,  J^{\#} \xi]  + [X, \xi] \\
        & = & [X,  J^{\#} \xi] +  ad^{*}_X \xi -  ad^{*}_{\xi}X  \\
                     & = & \{ [X,  J^{\#} \xi] -  J^{\#} ( ad^{*}_X \xi) \}
                     \,+ \,  \{ J^{\#} ( ad^{*}_X \xi) +  ad^{*}_X \xi \}.
\ee
Here we used the fact that  $ad^{*}_{\xi}X =0$, which can be 
easily verified directly.
It is easy to see that $ ad^{*}_X \xi  \in \f{h}^{\perp}$ so
that $J^{\#} ( ad^{*}_X \xi) +  ad^{*}_X \xi \in W_0$.
 Thus,  $[ X,  J^{\#} \xi + \xi] \in W_0$ iff
$$
 [X,  J^{\#} \xi] - J^{\#} ( ad^{*}_X \xi) =
  (ad_X \circ J^{\#}  - J^{\#} \circ ad^{*}_X ) \xi =
[ X, J]^{\#} \xi \, \in \f{h}.$$
Equivalently,
\begin{equation}
\label{J2}
[ X, \, J ] \equiv 0 \,(mod \,\f{h}) , \ \ \forall X \in \f{h},
\end{equation}
i.e., $J$ is \adh \,(mod $\f{h} $). Notice that,
as an element of $\f{g} \wedge \f{g}$, $J$ can always
 be written  as:
$$ 
J = \sum_{\alpha, \beta \in \Delta} J_{\alpha, \beta }
E_{\alpha}\wedge E_{\beta} \, + J_1
$$
where $ J_1 \equiv 0$ (mod \, $\f{h})$.
In fact, one can  always take $J_1 = 0$, which will  not
affect   the  Lagrangian subalgebra $W_0$. Moreover, it
follows from the equation:
$$    [ h, E_{\alpha}\wedge E_{\beta} ] =
     [ h, E_{\alpha}] \wedge E_{\beta}
    + E_{\alpha}\wedge [ h, E_{\beta} ] = <\al + \beta , h>
E_{\alpha}\wedge E_{\beta}  ,
\ \ \ \ \forall h \in \f{h},
$$
that  $ J_{\alpha, \beta } = 0$ whenever  $ \al + \beta \neq 0$.  By
denoting  $ J_{\alpha, -\al }$ by $J_{\alpha}$,   we  can write
\begin{equation}
\label{J3}
 J = \sum_{\alpha \in \Delta_{+}} J_{\alpha}
E_{\alpha}\wedge E_{-\alpha}, 
\end{equation}
which is in fact  \adh. \,
Thus, under the standard  identification
that $ \f{g}\oplus  \f{g^*}\cong \f{d} (= \f{g} \oplus \f{g})$, 
 $W_0$ is of the form (comparing  with Equation (\ref{B1}) in the
 last section):
\begin{equation}
\label{eq:JJ}
W_0= \{ (h, h), \, (X ,\varphi X ), \, (Y_- , Y_+) \, ~
|\, h \in \f{h},  \,X \in {\frak n^{\circ}_{\pm}}, \,
Y_{\pm} \in{\frak k_{\pm}} \},
\end{equation}
where 
$\f{k_{\pm}} 
= span_{\complex} 
 \{ E_{\pm \alpha} \, | J_{\alpha} =0, \alpha \in \Delta_{+} \}$,
$\f{n^{\circ}_{\pm}} =  span_{\complex} \{ E_{\pm \alpha} \, |
\, J_{\alpha} \neq 0, \alpha \in \Delta_{+} \}$  in  analogous to
Equation (\ref{ker}), and
 $\va$ is the Cayley transformation
of $(J^{\#} + r_0^{\#})|_{\frak n^{\circ}_{\pm}}$: 
$\varphi(\lambda) = { {(J^{\#} + r_0^{\#})+1} \over {(J^{\#}  + r_0^{\#})-1}}$.
Using a   similar argument as  in the proof of 
Lemma \ref{le:S}, we can show  that
$\f{n^{\circ}_{\pm}}$ are indeed subalgebras of $\f{n_{\pm}}$
and $\f{k_{\pm}}$ are ideals of   $\f{n_{\pm}}$.
Consequently, they  correspond to 
 a subset  $S$ of the set of simple roots  according to Proposition \ref{kn}.

Finally, by using  the fact that the commutator of the elements
$(E_{\alpha},\varphi_{\alpha} E_{\alpha})$ is  still in $W_0$,
 one derives the following relations:
$$
 \va_{-\alpha}=  \va_{\alpha}^{-1},  \ \ \ \ \va_{\alpha+\beta}
   = \va_{\alpha} \va_{\beta}, \ \ \ \ 
\forall \alpha, \beta \in \pm [S] \ \ \mbox{ such that }
 \alpha +\beta \in \Delta.
$$
This implies  that
$\va  = Ad_{e^{\lam_0}}$ for some $\lam_{0} \in \f{h}$.
This concludes the proof. \qed

In the sequel, we use  $l(S, \lam_0)$  to denote the Lagrangian
subalgebra $ W_0$ corresponding to the pair $(S, \lam_0)$.
Combining the above proposition and Corollary \ref{cor:ef} leads to:

\begin{thm}
\label{classify}
There is a one-one correspondence among the following objects:
\begin{enumerate}
\item dynamical $r$-matrices with zero gauge term,
\item  pairs $( S, \lam_0 )$, where $S$ is
 a subset  of the simple roots and $\lam_0 \in \f{h^*}$, and
\item Lagrangian subalgebras  $ W_0 \subset \f{g} \oplus \f{g}$ such that
 $W_0 \cap \f{g} = \f{h}$.
\end{enumerate}
\end{thm}

This  theorem  establishes the correspondence between
Lagrangian subalgebras  of $ \f{g} \oplus \f{g}$ and dynamical $r$-matrices
in a rather indirect manner, namely through the pair
$(S, \lam_0 )$.
Next we will illuminate a direct  connection geometrically.

Consider a   subbundle $\D$ of $A$, where
$U\subset TU$ is identified with the zero section and
$\frakh\subset \frakg$.
Given a function  $\T:  U \lon  \wedge^{2} \f{g} $,  being
considered  as a section in $ \Gamma ( \wedge^{2} A)$,
 the characteristic pair $(\D , \T )$ 
 defines a maximal isotropic subbundle $W$ of $A \oplus A^*$:
\begin{equation}
\label{eq:W}
  W = \{ \,X + \T^{\#} \xi + \, \xi \,| \,  
 X \in \D , \xi \in (\D )^{\perp} \},
\end{equation}
as given  by Equation (\ref{eq:L}). Then we have
\begin{pro}
\label{pro:W}
If $r (\lambda ) = \T (\lambda ) + r_0$ is a dynamical $r$-matrix, 
the subbundle $W$ corresponding to the characteristic pair $( \D ,\, \T )$ 
 is a Dirac structure of the  Lie bialgebroid $(A, A^*, \, r_0)$.
\end{pro}
\pf
 It suffices to   check the three conditions in Theorem \ref{th:pair}. First,
 it is obvious  that $\D \subset A$ is  a Lie subalgebroid.
Second, we have
\be
d_{*}\T  +\half [\T, \T]
         &=&[ r_0 , \T ] + \half[\T, \T] \\
        &=& - \sum h_i \wedge \plamm{\T} 
        \equiv 0 \, \ (mod \, U \times \f{h}),
\ee
according to Equation  (\ref{eq:t}).

Third, $\forall \xi, \eta \in  \Gamma((U \times \f{h})^{\perp})$
 and $h \in \f{h}$,

\be
< L_{\T^{\#}\xi} \eta, h > &=& < \eta \,, \,[ h , \T^{\#}\xi ] > \\
                     &=& < \eta , \T^{\#}(L_h \xi) > \\
                     &=& < [ h, \T^{\#} \eta ] , \xi > - L_h < \T^{\#} \eta , \xi > \\
                     &=& < L_{\T^{\#}\eta} \xi , h> + < d<\T^{\#}\xi , \eta > , h>, 
\ee
where in the second equality we used the fact that $\T$  is $\frakh$-invariant.
It thus follows  that
$$
<[\xi , \eta ] _{\T} , h>
=< L_{\T^{\#}\xi}\eta -L_{\T^{\#}\eta }\xi - d \langle \T^{\#}\xi , \eta \rangle , h >  = 0,
\ \ \ \ \forall h \in \f{h}.
$$
That is, $\gm  (U \times \f{h})^{\perp}$ is closed under
   $[\cdot , \cdot ]_{\T}$.
On the other hand,  it is  well known that
$\f{h}^{\perp}$
is an ideal of the dual Lie algebra $ \f{g}^* $,
 since $\f{h}\subset \frakg$ is a  Cartan subalgebra.
 This  means that
$(\D )^{\perp}$ is a Lie  subalgebroid of $A^*$. 
Thus,  
$\Gamma( U \times \f{h})^{\perp}$ is closed under
the bracket $[\xi , \eta ]+[\xi , \eta ]_{\T}$.
Consequently, the conclusion follows.  \qed

It is well known  that a Lie bialgebra integrates to
a Poisson group. Similarly, the global object corresponding to
a Lie bialgebroid is a Poisson groupoid \cite{MX94} \cite{MX98}.
For the Lie bialgebroid $(A, A^{*}, r_{0})$,
its Poisson groupoid is very simple to describe.
As a groupoid, it is simply the product of the pair groupoid $U\times U$ with
the Lie group $G$, where $G$ is a Lie group
with Lie algebra $\frakg$.
The Poisson structure is the product
of the zero Poisson structure on $U\times U$ with
the Poisson group structure on  $G$ defined by the $r$-matrix $r_{0}$.
According to Theorem 8.6 in \cite{LWX1},  the Dirac structure $W$ 
 corresponds to  a Poisson homogeneous  space $Q$
 of this  Poisson groupoid.
As a manifold, 
$$   Q = ( U \times U \times G)/ (U \times H) \cong U \times G/H,
$$
where $H \subset G$ is a closed subgroup with Lie algebra $\f{h}$.
It is not difficult to see that
for each fixed $\lambda \in U $,
$\{ \lambda \} \times G/H $ is a Poisson submanifold,  whereas
  the Poisson tensor is
$$
\pi_{Q }(\lambda )= p_{*} ( r_0^L -r_0^R + \T^L  (\lam)) =
p_{*} ( r^L(\lambda) -r_0^R).
$$ Here $p: G\lon G/H$ is
the projection, 
 $r^{L}(\lambda )$ refers
to the bivector field on $G$ obtained by the left translation
of $r(\lambda )\in \wedge^{2}\frakg $,
and $r^{R}_{0}$ refers
to the bivector field on $G$ obtained by the right  translation
of $r_0\in \wedge^{2}\frakg $.
 It is  simple  to see that
 $( G/H, \pi_Q( \lambda ))$ is a Poisson homogeneous  $G$-space.
Thus in this way
 we obtain a family of Poisson homogeneous $G$-spaces  parameterized
by $\lambda \in U$. It is not surprising that this  is the family
of Poisson homogeneous spaces studied by Lu \cite{Lu}. 

The corresponding family of  Lagrangian subalgebras   (or Dirac structures)
of the  Lie bialgebra $(\f{g}, \f{g^*} , r_0 )$  is just the fibers of $W$:
\begin{equation}
\label{eq:W1}
  W(\lam ) = \{ \,X + \T^{\#}(\lam )\xi + \, \xi \,| \,  
 X \in \f{h} , \xi \in \f{h}^{\perp} \} .
\end{equation}
In other words, $W(\lambda )$ 
corresponds to the  characteristic pair $(\f{h}, \, \T(\lam ))$.  
In fact, it is easy to see that
\begin{equation}
\label{ll}
W(\lam ) = l( S, \lam + \lam_0 ),
\end{equation}
where $(S, \lam_{0})$ is  the pair corresponding to
the dynamical $r$-matrix $r( \lambda )$ as  in Theorem \ref{classify}.
We now summarize the above discussion in the following
two corollaries.
\begin{cor}
\label{lagrange}
The following two statements are equivalent:
\begin{enumerate}
\item  The subbundle $W$ defined by the  characteristic pair
 $(U \times \f{h}, \T)$ is  a Dirac structure of the
 Lie bialgebroid $(A, A^*, r_0)$.
\item  For any fixed $\lam \in U$,  $W(\lam )$ defined by
 the  characteristic pair  $( \f{h}, \T (\lambda))$ is a Dirac structure for the
Lie bialgebra $(\f{g}, \f{g^*}, r_0)$.
\end{enumerate}
\end{cor}

\begin{cor} \cite{Lu}
\label{cor:W}
A dynamical $r$-matrix $r(\lambda ) $
 defines  a family of  Dirac structures $W(\lambda )$ of  the
Lie bialgebra $(\f{g}, \f{g^*}, r_0)$, which in turn
 corresponds to a family of  Poisson homogeneous  $G$ -spaces
 $( G/H, \pi_Q( \lambda ))$.
\end{cor}

Such a  family of  Lagrangian subalgebras is   said 
 to be governed by a dynamical $r$-matrix.
From Corollary \ref{lagrange}, we see
that the inverse of Proposition \ref{pro:W} is
not necessary  true,
because  $W$ being a Dirac structure is only a fiberwise property
 without involving  any dynamical  relation.
In fact, given a family of Lagrangian subalgebras $W(\lambda ), \ \forall
\lambda \in U$,  we may write
 $W(\lam) = l(S_{\lambda } ,\, \psi (\lam))$ for $\psi (\lam)\in \frakh^*$.
From Equation (\ref{ll}),
it follows that  $W(\lambda )$ is governed by a dynamical
$r$-matrix iff $S_{\lambda}$ is independent of 
$\lambda $ and
$\psi : \f{h} \lon \f{h}$ 
is a linear translation: $ \psi (\lam) = \lam + \lam_0$
 for some $\lambda_{0}\in \f{h}$.
Consequently, we have
\begin{cor}
Let $ \mu \in U$ be any fixed point, $W_0$
a Lagrangian subalgebra  of $\frakg \oplus \frakg$
 such that  $W_0 \cap \f{g} =\f{h} $. Then $W_0$ extends uniquely
 to a family of Lagrangian subalgebras $W(\lambda )$
such that  $W(\mu )=W_0$, which is governed by a  dynamical
$r$-matrix.
\end{cor}
\pf
Assume that  $W_0 =l(S, \lam_0)$.
Consider the pair $( S,  \lam_0 - \mu)$. This  corresponds  to
a dynamical $r$- matrix $r(\lambda )$  according to
 Theorem \ref{classify}. Let $W(\lambda )$ be its corresponding
 family of Lagrangian subalgebras. Then
 $ W(\lam) = l( S, \lam -\mu + \lam_0 )$.
 Thus   $W(\mu) = l(S, \lam_0) =W_0 $.   Moreover, 
it is clear that such an extension is unique.
\qed


\end{document}